\documentclass[12pt]{article}
\usepackage{amssymb, amsmath, amsthm, mathrsfs}
\usepackage[T2A]{fontenc}
\usepackage[english,russian]{babel}
\begin{document}
\small
\vskip3mm
\begin{center}
{\bf CARLEMAN-VEKUA EQUATION WITH A SINGULAR POINT} \vskip 3mm
{\bf A. Tungatarov }\\

{\it   tun-mat@list.ru}
\end{center}

{\it \textbf{Abstract.} In this article unconditional solvability
of the Carleman-Vekua equation with a singular point is proved,
the Riemann-Hilbert problem is solved, integral representations of
solutions, the structures  of their zeros and poles are received.


\textbf{AMS. Subject classifications: }30G20, 35J70, 31A10,
35C15.}
 \vskip 3mm

{\bf Introduction.}

Let $G$  be a bounded domain of the complex plane E with  boundary
$\Gamma\in C^{1,\alpha}$, $0<\alpha<1$,
   with the inner point $z=a$.

Let $S(G)$ be  the set of measurable, essentially bounded
functions $f(z)$  in $G$ with the norm
$$||f||_{S(G)}=ess\sup_{z\in G}|f(z)|=\lim_{p\rightarrow \infty}||f||_{L_p(G)}.$$
Now the spaces used below are defined:

$S_\nu (G,a)$ is the class of functions $f(z)$, for which
$f(z)|z-a|^\nu\in S(G)$. The norm in $S_{\nu}(G,a)$  is defined by
the formula
$$||f||_{S_\nu(G,a)}=ess\sup_{z\in G}(|z-a|^\nu |f(z)|) ,$$
where $\nu$ is a real number.

  $C_\nu(\overline{G},a)$ is the class of functions $f(z)$, for which $f(z)|z-a|^\nu\in C(\overline{G})$.
  The norm in  $C_{\nu}(\overline{G},a)$ is defined by the formula
  $$||f||_{C_\nu(\overline{G},a)}=\max_{z\in G}(|z-a|^\nu |f(z)|).$$

 $U_0(G)$  is the class of  holomorphic  functions in $G$.

 $W_p^1(G)$, $p>1$   is the  Sobolev space, see [1].

Let us consider  the equation
 \begin{equation}\label{eq1}
 \partial_{\bar{z}}V+A(z)V+B(z)\overline{V}=F(z),
 \end{equation}

in  $G$,  where
$\partial_{\bar{z}}=\displaystyle\frac{1}{2}\bigg(\frac{\partial}{\partial
x}+i\frac{\partial}{\partial y}\bigg)$; $A(z),B(z)\in S_1(G,a)$;
$F(z)\in S_\beta(G,a)$, \,$\max \{0, 1 - 8\mu \}< \beta <
\frac{2}{q}.$     Here $\mu =
\|A\|_{S_{1}(G,a)}+\|B\|_{S_{1}(G,a)}; \,\,$ $2<q<\frac{2}{1-8
\mu}$ if $\mu <\frac{1}{8}$  and  $q>2$ if $\mu \geq \frac{1}{8}.$

For $a=0$  we receive the equation

 \begin{equation}\label{eq2}
\partial_{\bar{z}}V+\frac{A_0(z)}{|z|}V+\frac{B_0(z)}{|z|}\overline{V} =F(z),
 \end{equation}

where  $A_0(z), B_0(z)\in S(G)$; $F(z)\in S_\beta (G,\, 0).$

It is obvious, that $F(z)\in L_q(G), \,\, q>2$.

The equation (2) with $F(z)\equiv0$ arise in the theory of infinitesimal deformations of surfaces
of positive curvature with a flat point, see [2, 3]. There it is required to prove the existence of
continuous solutions of equation (2) in the neighborhood of a singular point $z=0$. In this connection
the equation (2) is studied in many works of L.G. Mikhailov, Z.D Usmanov, A. Tungatarov, see [4-7], etc.

In all results usually a sufficient smallness of the coefficients  $A_0(z)$ and $B_0(z)$ or the
 smallness of the domain $G$ is supposed. In the present article we prove existence of solutions
  $V(z)$ of the equation (2), satisfying to the condition $V(0)=0$, without any conditions on the
   smallness of the coefficients or on the  smallness of the domain $G$. Integral representations
   of  solutions of equation  (1) in the class

 \begin{equation}\label{eq3}
W_q^1(G)\cap C_{\beta-1}(\overline{G},a), 0<\beta<\frac{2}{q}, q>2
 \end{equation}
 are also received.  The Riemann  -- Hilbert problem for the equation (2) is solved in the class (3)
 without any smallness conditions on the coefficients assumed in [4-7]. The structures of the zeros
 and poles of the  solutions of the equation (1)  in the class (3) are investigated.

\begin{center}
{\bf \S 1. Representations of solutions and solvability of the
equation (1) }
\end{center}

The solutions of the equation (1) are looked for in the class  (3). As in  [1] it is possible to prove, that the  equation  (1) is equivalent to the equation
  \begin{equation}\label{eq4}
V(z)+(P_GV)(z)=(T_GF)(z)+\Phi(z), \,\ z\in G,
 \end{equation}
 where   $$(T_Gf)(z)=-\frac{1}{\pi}\iint\limits_G{\frac{f(\zeta)}{\zeta-z}}dG_\zeta, \,\,\, (dG_\zeta=d\xi d\eta , \,\,\, \zeta=\xi+i\eta), $$
 $$\Phi(z)\in U_0(G), \,\ (P_G f)(z)=(T_Gf^{*})(z), \,\ f^{*}=A(z)f+B(z)\bar{f}.$$

 As the function  $V(z)$ belongs to the class  (3), then $V(a)=0$. Therefore  from  (4) we receive

 \begin{equation}\label{eq5}
 (P_GV)(a)=(T_GF)(a)+\Phi(a).
  \end{equation}
From  (4) and (5) it follows
\begin{equation}\label{eq6}
V(z)+(P_{G,a}V)(z)=(T_{G,a}F)(z)+(z-a)\Phi(z),\,\ z\in G,
\end{equation}
 where
$$(P_{G,a}f)(z)=(P_Gf)(z)-(P_Gf)(a),$$ $$(T_{G,a}f)(z)=(T_Gf)(z)-(T_Gf)(a),\,\ \Phi(z)\in U_0(G).$$

Thus, any solution of  equation (1) from the class (3) satisfies the equation (6). A solution of the equation (6) is looked for in the class $C_{\beta-1}(\overline{G},a)$. Let us  show, that  any solution of the equation (6) from the class  $C_{\beta-1}(\overline{G},a)$ belongs to the class  (3) and almost everywhere in $G$   satisfies the equation  (1). Hereafter $M$  denotes positive constants, not depending on the factor involved.
\vskip 3mm
{\bf Lemma 1.} {\it The operator  $(T_{G,a}f)(z)$  maps the class  $S_\beta(G,a)$, $\displaystyle 0<\beta<\frac{2}{q}$, $q>2$   into the class  $C_{\beta-1}(\overline{G},a)\bigcap C^\alpha(\overline{G})$, $\alpha=1-\displaystyle\frac{2}{q}$, $\displaystyle 0<\beta<\frac{2}{q}$, $q>2$  and moreover, the estimates

 \begin{equation}\label{eq7}
 ||(T_{G,a}f)(z)||_{G_{\beta-1}(\overline{G},a)}\leq M_\beta(G)||f||_{S_\beta(G, \, a)},
 \end{equation}

  \begin{equation}\label{eq8}
 ||(T_{G,a}f)(z)||_{C^\alpha(\overline{G})}\leq M ||f||_{L_q(G)}
 \end{equation}
are hold, where
 $$M_\beta(G)=\sup_{z\in G}\frac{|z-a|^\beta}{\pi}\iint\limits_{G}\frac{dG_{\zeta}}{|\zeta-a|^{1+\beta}|\zeta-z|}.$$}

 \vskip 3mm
{\bf Proof.}  Let  $f(z)\in S_\beta (G,a)$, $\displaystyle 0<\beta<\frac{2}{q}$, $q>2$. Using the Hadamard  inequality, we have
 $$|(T_{G,a}f)(z)|\leq M_{\beta}(G)||f||_{S_{\beta}(G,a)} |z-a|^{1-\beta}.$$

From here  the inequality  (7) follows. As  $f(z)\in L_q(G)$, $q>2$, the inequality  (8) follows from the estimation  (6.8) of the work  [1, ch.1, \S 6]. The lemma is proved.
\vskip 3mm
{\bf Lemma 2.} {\it The  operator $(P_{G,a}f)(z)$  maps the space $C_{\beta-1}(\overline{G},a)$   into the space  $C_{\beta-1}(\overline{G},a)\bigcap C^\alpha(\overline{G})$, $\displaystyle\alpha=1-\frac{2}{q}$, $0<\beta<\frac{2}{q}$, $q>2$, and moreover,  the estimates
  \begin{equation}\label{eq9}
 ||(P_{G,a}f)(z)||_{G_{\beta-1}(\overline{G},a)}\leq \mu M_\beta(G)||f||_{C_{\beta-1}(\overline{G},a)},
 \end{equation}

  \begin{equation}\label{eq10}
 ||(P_{G,a}f)(z)||_{G^\alpha(\overline{G})}\leq \mu M ||f||_{C_{\beta-1}(\overline{G},a)},
 \end{equation}
  \begin{equation}\label{eq11}
 ||(P_{G,a}f_1)(z)-(P_{G,a}f_2)||_{G_{\beta-1}(\overline{G},a)}\leq \mu M_\beta (G)||f_1-f_2||_{C_{\beta-1}(G,a)},
 \end{equation}
 are true,  where   $\mu=||A||_{S_1(G,a)}+||B||_{S_1(G,a)}$;  $f_1, \, f_2$ are arbitrary functions from the class  $C_{\beta-1}(\overline{G},a)$.}
\vskip 3mm

{\bf Proof.}   Let  $f(z)\in C_{\beta-1}(\overline{G},a)$. Then $
f^{*}(z) \in S_{\beta}(G,a)\subset_{\rightarrow} L_{q}(G),
2<q<\frac{2}{\beta} $ and the estimation
$$ \| f^{*}\|_{S_{\beta}(G,a)}\leq \mu \| f \|_{C_{\beta-1}(G,a)}
$$
is hold. Therefore, from (7) the inequality (9) follows. The inequality (11) can be proved similarly. As
$$ \| f^{*} \|_{L_{q}(G)}\leq M\mu \|f \|_{C_{\beta-1}(G,a)}
$$
inequality (10) follows from the estimation (6.8) of the work [1, ch.1, \S 6]. The lemma is proved.

From the form of  the  equation  (6) and lemmas 1 and 2 by virtue of  the results in   [1], it follows, that any solution of this equation from the class  $C_{\beta-1}(\overline{G},a)$ belongs to the class $W_q^1(G)$, $0<\beta<\frac{2}{q}$, $q>2$, and satisfies the equation (1) almost everywhere in $G$. Thus, the next result is true.

\vskip 3mm

{\bf Theorem 1.} {\it  Any solution of equation  (1) from the  class (3) satisfies the equation  (6). And vice  versa,  if \  $\Phi(z)\in U_0(G)$,  then any solution  of the equation (6) from the class $C_{\beta-1}(\overline{G},a)$ belongs to the class $W_q^1(G)$, $\displaystyle 0<\beta<\frac{2}{q}$, $q>2$,   and satisfies the equation (1) almost everywhere in  $G$ .}

Let
$$(K_\Gamma f)(z)=\frac{1}{2\pi i}\int\limits_{\Gamma}\frac{f(t}{t-z}dt, \,\ (K_{\Gamma,a }f)(z)=(K_\Gamma f)(z)-(K_\Gamma f)(a).$$

 Applying the operator $(K_\Gamma f)(z)$   for  $z\in G$  to both parts of equality  (6), we receive

 \begin{equation}\label{eq12}
 (K_\Gamma V)(z)+(P_{G}V)(a)=(T_GF)(a)+(z-a)\Phi(z).
 \end{equation}
From here for  $z=a$  we have

\begin{equation}\label{eq13} (K_\Gamma V)(a)+(P_G
V)(a)=(T_GF)(a).
 \end{equation}
From (12) and (13) it follows
 $$(z-a)\Phi(z)=(K_{\Gamma,a}V)(z).$$
Hence, the integral representation of the first type for the solutions of equation (1)  from the class (3) has the form
 $$V(z)=-(P_{G,a}V)(z)+(K_{\Gamma,a}V)(z)+(T_{G,a}F)(z), \,\ z\in G . $$

\vskip 3mm
{\bf Theorem 2.} {\it The equation (6) is solvable  in the class $C_{\beta-1}(\overline{G},a)$ for any right-hand side from the same class. }

{\bf Proof. }
 As $0<\beta<1$, then from $V(z)\in C_{\beta-1}(\overline{G},a)$ it follows, that  $V(z)\in C(\overline{G})$ and
 \,\,\,\,\,\,
$V(a)=0$. Therefore, from lemma 2 and Arzela's theorem it follows,
that $P_{G,a}$  is  a completely  continuous operator from
$C_{\beta-1}(\overline{G},a)$  into
$C_{\beta-1}(\overline{G},a)\bigcap C^\alpha(\overline{G})$.
 Hence, for the proof of the solvability of equation (6) in class  $C_{\beta-1}(\overline{G},a)$  it is enough to show, that
 the corresponding homogeneous equation
\begin{equation}\label{eq14}
V(z)+(P_{G,a}V)(z)=0, \,\ z\in G
 \end{equation}
has only the trivial  solution in $C_{\beta-1}(\overline{G},a).$

Let us prove now, that the homogeneous equation (14) has only the trivial solution in class $C_{\beta-1}(\overline{G},a)$.
 The proof is performed by contradiction.

Assume that equation (14) has a not-trivial solution  $V(z)$ in the class $C_{\beta-1}(\overline{G},a)$.

Applying the operator $ \partial_{\overline{z}} $ to both sides of
the equation (14) we have
$$\partial_{\overline{z}}V+A(z)V+B(z)\overline{V}=0. $$
Hence, the representation for solution
\begin{equation}\label{eq15}
V(z)=\Phi(z) \exp(-\omega(z))
 \end{equation}
 holds, see [1,4], where $\Phi(z)\in U_{0}(G), \, \omega(z)=(T_{G}V^{\wedge})(z),$ $V^{\wedge}=\frac{V^{*}}{V} $.

From $(15)$  it follows
  \begin{equation}\label{eq16}
\Phi(z)=V(z) \exp\omega(z)
 \end{equation}

Substituted $(15)$  into $(14)$ gives

 $$ \Phi(z)=-\exp(-\omega(z))\cdot (P_{G,a}V)(z), \,\,  where \,\,\, V=\Phi(z) \exp(-\omega(z)). $$

As  \,\, $ \exp(-\omega(z))\cdot (P_{G,a}V)(z)\in U_{0}(E
\setminus \overline{G}),$ then from  the  last  equality  it
follows that the  function $\Phi(z) $ analytically  single extends
to the  entire complex plane   $E$.

Let $0<\varepsilon<1$ and $G_{0}=\{z:|z-a|<\varepsilon\}\subset G.$ Using the Hadamard inequality, we have

$$|\omega(z)|\leq\frac{\mu}{\pi}(\iint\limits_{G \setminus G_{0}}\frac{dG_\zeta}{|\zeta-a||\zeta-z|}
-\iint\limits_{G_{0}}\frac{dG_\zeta}{|\zeta-a||\zeta-z|})\leq
\mu(M(G \setminus G_{0})+8\ln\frac{1}{|z-a|}),  $$ where $M(G
\setminus G_{0})>0 $ is a constant, depending only on $G \setminus
G_{0}$. From here it follows

$$ M|z-a|^{8\mu} \leq |e^{\omega(z)}| \leq \frac{M}{|z-a|^{8\mu}},
z\in G $$

\,\,\,Hence by virtue of (16) it follows that $ \Phi(z)\in
C_{8\mu+\beta-1}(\overline{G},a). $

As  $ \beta>1-8\mu $  if  $ 1-8\mu>0 $  and  $ \beta>0 $
  if  $ 1-8\mu\leq 0 $  then from here  it follows that $ 8\mu+\beta-1>0 $
and $ \Phi(\infty)=0 $.

Therefore, the function $ \Phi(z)$ analytically single extends to
the entire complex plane and is equal to zero at point $z=\infty
$. Then by Liouville theorem $ \Phi(z)\equiv0 $. Therefore, $
V(z)\equiv 0 $ in $G$.  The theorem is proved.

From this theorem by virtue of theorem 1 the next result it
follows.

\vskip 3mm {\bf Theorem 3.} {\it The equation  (1) is solvable in
the class (3).} \vskip 3mm

\begin{center}
{\bf \S 2. Zeros and poles of solutions of  equation  (1)}
\end{center}

Let $k$  be an integer number. Let us consider the equation  (1)
in $G$, where $F(z)\in S_{\beta-k}(G,a)$, $0<\beta<1$; $A(z),
B(z)\in S_1(G,a)$.

The solution of equation (1) from the class
 \begin{equation}\label{eq17}
W_q^1(G)\bigcap S_{\beta-k-1}(G,a) \end{equation} is looked in the
form
  \begin{equation}\label{eq18}
V(z)=(z-a)^kW(z),
\end{equation}
where  $W(z)$ is a new unknown function from the class (3).

Substituted (18) into (1) we get
 \begin{equation}\label{eq19}
\partial_{\bar{z}}W+A(z)W+B_k(z)\overline{W}=F_k(z),
\end{equation}
where
$$B_k(z)=B(z)\exp(-2ik\varphi), F_k(z)=(z-a)^{-k}F(z), \, \varphi=\arg (z-a).$$

It is obvious, that  $B_k(z)\in S_1(G,\, a)$, $F_k(z)\in
S_\beta(G, \, a)$.
 Therefore, by virtue of the results of \S1 the equation (19) has  solutions from the class (3). They can be found from the equation
  \begin{equation}\label{eq20}
W(z)+(P^{\wedge}_{G,\, a}W)(z)=(T_{G,\, a}F_k)(z)+(z-a)\Phi(z),
\end{equation}
where
$$(P_{G,\, a}^{\wedge}f)(z)=(T_{G,a}f^{*}_k)(z),\,\ f_{k}^{*}(z)=A(z) f+B_k(z) \bar{f}.$$

Thus the next theorem is proved.

\vskip 3mm {\bf Theorem 4. } {\it The equation (1), where $F(z)\in
S_{\beta-k}(G,a)$, $\displaystyle 0<\beta<\frac{2}{q}$, $q>2$, $k$
is an integer number, $A(z), B(z)\in S_1(G,a)$,  has solutions
from the class (17), which can be found by formulas (18), (20).}
\vskip 3mm

\begin{center}
{\bf \S 3. Generalized Riemann-Hilbert problem for equation (1)}
\end{center}

Let  $R>0$, $G=\{z:|z|<R\}$, $\Gamma=\{t:|t|=R\}$.  Let us consider  equation  (1) in $G$, where  $A(z), B(z)\in S_1(G,0)$, $F(z)\in S_\beta(G,0)$, $\displaystyle 0<\beta<\frac{2}{q}$, $q>2$.

We look for solutions of equation (1)in the class
\begin{equation}\label{eq21}
W_q^1(G)\bigcap S_{\beta-1}(G,0), \,\ 0<\beta<\frac{2}{q}, \,\
q>2.
\end{equation}

Let us consider in $G $ the  generalized  Riemann-Hilbert problem in the canonical form, see [1]. The more general cases can be reduced to this form, see  [1].

{\bf Problem R--H}. {\it It is necessary to find a solution of
equation (1) in the class (21), satisfying the boundary condition
\begin{equation}\label{eq22}
Re[t^{-m}V(t)]=g(t), \,\ t\in \Gamma,
\end{equation}
where $m$ is an integer number,  $g(t)\in C^{\alpha}(\Gamma)$,
$0<\beta<\frac{2}{q},\,\ q>2$. }

 {\bf $1^{0}$.}  Let   $m\geq 1$. To solve the R--H problem the equation (6) is used with $a=0$:
\begin{equation}\label{eq23}
V+(P_{G,0}V)(z)=(T_{G,0}F)(z)+z\Phi(z), \,\ z\in G,
\end{equation}
where    $\displaystyle\Phi(z)\in U_0(G)\bigcap C^{\alpha}(\overline{G}), \alpha=1-\frac{2}{q}$.

Following  [1, ch.4, \S7],  the function $\Phi(z)$ in  equation
(23) is represented in the  form
\begin{equation}\label{eq24}
\Phi(z)=\Phi_0(z)-(P_mV)(z)+(Q_mF)(z),
\end{equation}
where
 $$\Phi_0(z)\in U_0(G)\bigcap C^{\alpha}(\overline{G}), \,\ (P_mV)(z)=(Q_mV^{*})(z), \,\ V^{*}=A(z)V+B(z)\overline{V}, $$
 $$(Q_mf)=-\frac{z^{2m}}{\pi R^{2m}}\iint\limits_{G}\frac{\overline{f(\zeta)}dG_\zeta}{R^2-\bar{\zeta}z}-\frac{z^{2m-1}}{\pi R^{2m}}\iint\limits_{G}\frac{\overline{f(\zeta)}}{\bar{\zeta}}dG_\zeta.$$

 Substituting the representation (23) for $V(z)$  in the boundary condition (22),
we have
$$
Re[t^{-m+1}\Phi_0(t)]=g(t).
$$

The general solution of this problem is given by the formula, see  [1, ch. 4, \S7],

\begin{equation}\label{eq25}
\Phi_0(z)=(D_{m-1}g)(z)+\Phi_{0m}(z),
\end{equation}
where
$$
(D_mg)(z)=\frac{z^m}{2\pi i}\int\limits_{\Gamma}g(t)\frac{t+z}{t-z}\frac{dt}{t},
$$
 \begin{equation}\label{eq26}
  \Phi_{0m}(z)=\sum\limits_{k=0}^{m-2}(\alpha_k(z^k-R^{2(k-m-1)}z^{2m-k-2})+i\beta_k(z^k-R^{2(k-m+1)} z^{2m-k-2}))+i\beta_m z^{m-1},
 \end{equation}
if $m\geq 2$ and $ \Phi_{0m}(z)=i\beta_m,\,\ if \,\  m=1. $

Here  $\alpha_k, \beta_k$, $k=0,...,m-2$; $\beta_m$    are arbitrary real numbers.

From formulas (23) - (25) it follows that
\begin{equation}\label{eq27}
V(z)+(P^{\wedge}_mV)(z)=(D_mg)(z)+(Q^{\wedge}_mF)(z)+z\Phi_{0m}(z),
\,\ z\in G,
\end{equation}

where
$$(P_m^{\wedge}V)(z)=(P_{G,0}V)(z)+z(P_mV)(z), \,\ (Q_m^{\wedge}V)(z)=(T_{G,0}F)(z)+z(Q_mF)(z).$$

Thus, with  $m\geq 1$ the R-H problem is reduced to the equivalent
equation (27). For any real numbers $a_k$  and  $\beta_k$,
$k=0,...,m-2$; $\beta_m$    the solution of equation (27) is the
solution of the R-H problem. Let us prove, that equation (27)  has
a solution in the class  $C_{\beta-1}(\overline{G},0)$. As in the
case of the operator   $(P_{G,a}V)(z)$ it is proved, that the
operator $(P_m^{\wedge}V)(z)$   is  completely continuous in the
class   $C_{\beta-1}(\overline{G},0).$

Therefore, our statement will be proved, if we show, that the homogeneous equation
 \begin{equation}\label{eq28}
V+(P_m^{\wedge}V)(z)=0,\,\ z\in G
\end{equation}
has no non-trivial solution in the class $C_{\beta-1}(G,0)$. From
(28) by  virtue of the Cauchy integral formula, see  [8], we get
$$(K_\Gamma V)(z)-(P_GV)(0)=-z(P_mV)(z), \,\ z\in G.$$

If we compare the coefficients of the series expansions with respect to $z$ of the left-and right-hand sides of the last equality, then we obtain
\begin{equation}\label{eq29}
\int\limits_\Gamma V(t)e^{-ik\theta}d\theta=0,\,\ k=0,...,2m-1;\,\
t=Re^{i\theta}.
\end{equation}

Moreover, any solution of equation  (29), can be represented as
\begin{equation}\label{eq30}
V(z)=\Phi(z) \exp\Omega(z),
\end{equation}
where
$$\Phi(z)\in U_0(G)\bigcap C^\alpha(\overline{G}), \,\ \Phi(0)=0,$$
  $$\Omega(z)=-\frac{1}{\pi}\iint\limits_G\left(\frac{V^{\wedge}(\zeta)}
  {\zeta-z}-\frac{z\overline{V^{\wedge}(\zeta)}}{R^2-\bar{\zeta}z}\right)dG_\zeta,\,\ V^{\wedge}(z)=\frac{V^{*}(z)}{V(z)}.$$

  But  $V(z)$ satisfies also the homogeneous boundary condition
\begin{equation}\label{eq31}
Re[t^{-m}V(t)]=0,\,\ t\in \Gamma.
\end{equation}
As ${Re[i\Omega(t)]=0}$, ${t\in \Gamma}$, then the boundary
condition  (31) by (30) is  represented as
 $$Re[t^{-m}\Phi(t)]=0.
 $$

The general solution of this problem with  $\Phi(0)=0$   has the form
 $$\Phi(z)=\sum\limits_{k=1}^{2m-1}c_kz^k,$$
where $c_k,\,\ k=1,...,2m-1$; are complex constants, satisfying
the conditions ${c_{2m-k}=-\bar{c}_k}$, ${k=1,...,m}$. Therefore,
by virtue of  (30) the solution of equation  (28) has the form
 $$V(z)=\sum\limits_{k=1}^{2m-1}c_kz^k\exp\Omega(z).$$
Inserting this into the equality  (29), we get
\begin{equation}\label{eq32}
\sum\limits_{k=1}^{2m-1}c_k
\int\limits_{\Gamma}t^kt^{-n}\exp\Omega(t)dt=0,\,\  n=1,...,2m-1.
\end{equation}

From here follows, that  $c_k=0$, $k=1,...,2m-1$; because the
determinant of the system (32) is different from zero as the Gram
determinant for the system of the linearly independent functions
$$t^k\exp\left(\frac{\Omega(t)}{2}\right),\,\ k=1,...,2m-1; \,\ \Omega(t)=\overline{\Omega(t)},\,\ t\in\Gamma.$$

This proves, that the homogeneous equation (28) has only the
trivial solution. Hence, the integral equation  (27) has a
solution in  $C_{\beta-1}(\overline{G},0)$   for any right -- hand
side belonging to $C_{\beta-1}(\overline{G},0)$.

Thus for  $m>0$  the non -- homogeneous R-H problem  is always
solvable a solution and the homogeneous R-H problem ($F\equiv0,\,\
g\equiv 0$)  has exactly $2m-1$  linearly independent solutions
over  the field of real numbers. The latter follows from the fact
that the homogeneous problem is equivalent to the integral
equation
$$V+(P_m^{\wedge}V)(z)=z\Phi_{0m}(z),$$
where $\Phi_{0m}(z)$  is defined by formula  (26),   the
right-hand side of which is  a linear combination of  ${2m-1}$
linearly independent  functions. Thus the next theorem follows.

\vskip 3mm {\bf Theorem 5.}  {\it  For  $m>0$ the non-homogeneous
R-H problem is  always solvable and the homogeneous R-H problem
has $2m-1$ linearly independent solutions over the field of real
numbers.} \vskip 3mm

{\bf {$2^0.$} } Let $m=0$. Following to $1^0$,  the function
$z\Phi(z)$ in the equation (23)  is represented in the form
  \begin{equation}\label{eq33}
z\Phi(z)=(D_0g)(z)+ic_0-z(P_0V)(z)+z(Q_0F)(z),
\end{equation}
where   $c_0$ is an arbitrary real number.

The  function $\Phi(z)$, given by formula (33),  belongs to the
class  $U_0(G)$, if the  equalities
  \begin{equation}\label{eq34}
(D_0g)(0)=0,
\end{equation}

 \begin{equation}\label{eq35}
c_0=0
\end{equation}
are hold.

Using (33)--(34), from (23) we obtain
 \begin{equation}\label{eq36}
V(z)+(P^{\wedge}V)(z)=(D_0g)(z)+(Q^{\wedge}F)(z), \,\ z\in G,
\end{equation}
where
$$(P^{\wedge}V)(z)=(P_{G,0}V)(z)+z(P_0V)(z),\,\ (Q^{\wedge}V)(z)=(T_{G,0}F)(z)+z(Q_0F)(z).$$

If equality (34) holds, then any solution of equation  (36)
belonging to $C_{\beta-1}(\overline{G},0)$ satisfies the boundary
condition (22) when $m=0$. Solvability of equation  (36) in the
class $C_{\beta-1}(\overline{G},0)$ can be proved similarly to the
proof of the  solvability  of  the equation  (27). Hence, under
the  condition (34) the R-H problem is solvable. Thus, we have the
next  result.

\vskip 3mm {\bf Theorem 6.} {\it For   the solvability of the  R-H
problem  in the case $m=0$ it is necessary and sufficient that the
single condition (34) is satisfied. If the condition (34) holds,
then the solution of the problem can be found from equation (36).
} \vskip 3mm

{\bf $3^0.$} Let $m<0$. As the solutions of the R-H problem are
sought from the class (21), in this case formula (27) is not
suitable. Therefore, let us introduce the function
$W(z)=z^{k+1}V(z)$ into consideration, where $k=-m$. The function
$W(z)$  satisfies the equation
$$\partial_{\bar{z}}W+A(z)W
+B_{-k-1}(z)\overline{W}=F_{-k-1}W(z),$$
where
 $$B_{-k-1}(z)=\exp(2(k+1)i\varphi) B(z), \,\ F_{-k-1}(z)=z^{k+1}F(z), \,\ \varphi=\arg z$$
and the boundary condition
$$Re[t^{-1}W(t)]=g(t),\,\ t\in \Gamma.$$
 It is obvious, that  $B_{-k-1}(z)\in S_1(G,0)$, $F_{-k-1}(z)\in S_{\beta}(G,0)$.

Therefore this problem corresponds to the one considered in $1^0$ for $m=1$. Hence, the function $W(z)$  satisfies the equation
 $$W(z)+(P_{1}^{\wedge}W)(z)=(D_1g)(z)+(Q_{1}^{\wedge}F_{-k-1})(z)+ic_0z,\,\ z\in G,$$
  where $c_0$   is an arbitrary real number,
$$(P_1^{\wedge}W)(z)=(T_{G,0}W^{*}_{-k-1})(z)+z(Q_1W^{*}_{-k-1})(z),\,\ W_{-k-1}^{*}=A(z)W+B_{-k-1}(z)\overline{W}.$$

Thus, when  $m<0$  the solution of the  R-H problem can be found from the equation
  \begin{equation}\label{eq37}
V(z)+(QV)(z)=\frac{z}{\pi
i}\int\limits_{\Gamma}\frac{g(t)dt}{t^{k+1}(t-z)}+(T_{G,0}F_{-k-1})(z)+\sum\limits_{j=1}^{k}a_jz^{-j},\,\
z\in G,
\end{equation}
where
$$(QV)(z)=(T_{G,0}(z^{k+1}V^{*}(z)))(z)-\frac{z}{\pi}\iint\limits_{\Gamma}\frac{\bar{\zeta}^{2k}\overline{V^{*}(\zeta)}}{1-\bar{\zeta}z}dG_{\zeta},$$

$$a_j=\frac{1}{\pi}\iint\limits_{G}\zeta^{j-1}V^{*}(\zeta)dG_{\zeta}+\frac{1}{\pi}\iint\limits_{G}\bar{\zeta}^{2k-j-1}\overline{V^{*}(\zeta)}dG_\zeta+\frac{1}{\pi i}\int\limits_\Gamma t^{j-k-1}g(t)dt,\,\ g=0,...,k-1;$$

$$a_k=\frac{1}{\pi}\iint\limits_{G}\zeta^{k-1}V^{*}(\zeta)dG_{\zeta}+\frac{1}{2\pi i}\iint\limits_{\Gamma}\frac{g(t)}{t}dt+ic_0,\,\ V^{*}=A(z)V+B(z)\overline{V}.$$

From (37) it follows  that for the  continuity of the function
$V(z)$ inside of $G$  it is necessary and sufficient that the
equalities
\begin{equation}\label{eq38}
a_j=0,\,\ j=0,...,k;
\end{equation}
are hold. The condition  (38) contain  $2k+2$ real equalities. One
of them, namely $Im  \ a_k=0$, is possible to be satisfied  by
means of a suitable choice of the constant $c_0$. Hence, there
remain  $2k+1$  conditions. Thus, for  $m<0$  the R-H  problem is
reduced  to the equation
\begin{equation}\label{eq39}
V(z)+(QV)(z)=\frac{z}{\pi
i}\int\limits_\Gamma\frac{g(t)dt}{t^{k+1}(t-z)}+(T_{G,0}F_{-k-1})(z),\,\
z\in G.
\end{equation}
The operator $Q$  is completely continuous in $C_{\beta-1}(G,0)$  and mapping this space into $C_{\beta-1}(\overline{G},0)\bigcap C^\alpha(\overline{G})$.

If in the equation
 $$V(z)+(QV)(z)=0$$
we replace $V(z)$ by $zW(z)$, then we obtain the equation (7.33)
from [1, ch.4, \S 7], which has only the trivial solution.
Therefore equation (39)
 is solvable in $C_{\beta-1}(G,0)$  for any right-hand side from the same class. Thus, we have the next result.

\vskip 3mm {\bf Theorem 7.} {\it For the solvability of the R-H
problem in the case  $m<0$ it is necessary and sufficient that the
${2|m| +1}$ real conditions  (38) are satisfied.}

\begin{center}
{\bf \S4. Riemann-Hilbert problem with an initial condition for equation (1)}
\end{center}

Let  $G=\{z:|z|<R \}$, $\Gamma=\{t:|t|=R\}$, $\nu>0$, $k=[\nu]$, $k\neq \nu$, $\beta=1-\nu+k$, $R>0$.
 Let us consider the equation (1) in $G$, where $A(z),\,\ B(z)\in S_1(G,0)$, $F(z)\in S_{1-\nu}(G,0)$.

Let us consider the  Riemann-Hilbert problem with an initial condition in the following form.

\vskip 3mm
{\bf Problem $(R-H)_0$}. {\it Find the solution of the equation (1) from the class
 $$C_{-\nu}(G,0)\bigcap W_q^1(G),\,\ 2<q<\frac{2}{\beta},$$
satisfying the boundary condition
 \begin{equation}\label{eq40}
 Re[t^{-n}V(t)]=g(t),\,\ t\in \Gamma,
 \end{equation}
where  $g(t)\in C^{\alpha}(\Gamma)$,
$\displaystyle\alpha=1-\frac{2}{q}$,   \ \   $n$ is an integer
number.}

The solution of the  $(R-H)_0$  problem   is looked for in the form
  \begin{equation}\label{eq41}
V(z)=z^{k}W(z),
\end{equation}
where  $W(z)$  is a new unknown function from the class
$$W_q^1(G)\bigcap S_{\beta-1}(G,0), \,\ \displaystyle 0<\beta<\frac{2}{q}, \,\ q>2.$$

\vskip 3mm {\bf Remark  1.} \ {\it If  \ $\nu<1$, the substitution
(41) is not required.} \vskip 3mm {\bf Remark 2.} \ {\it  If  \
$\nu>0$ is an integer number, then from
$V(z)=O(|z|^{\nu+\beta})$, $z\rightarrow 0$, $0<\beta<1$,  it
follows  that $V(z)=O(|z|^\nu)$, $z\rightarrow0$. Therefore, the
results which will  be obtained for   $[\nu]\neq \nu$, will also
hold for $[\nu]=\nu$.}

Substituting  (41) into (1) and (40), respectively, we obtain
\begin{equation}\label{eq42}
\partial_{\bar{z}}W+A(z)W+B_k(z)\overline{W}=F_k(z), \,\ z\in G,
\end{equation}
where
$$B_k(z)=B(z)\cdot \exp(-2ik\varphi), \,\ F_k(z)=z^{-k}F(z)$$
and
$$Re[t^{-m}W(t)]=g(t),\,\ t\in \Gamma,\,\ m=n-k.$$

It is obvious, that $B_k(z)\in S_1(G,0)$, $F_k(z)\in S_\beta(G,0)$, $0<\beta<1$. Hence, we obtain the
  Riemann -Hilbert problem solved in \S3 for the equation (42). Therefore,  from the results of \S3 the next result follows.

{\bf Theorem 8.} {\it 1) For  $n>[\nu]$ the problem $(R-H)_0$
 is always solvable.  The corresponding homogeneous problem has  $2(n-[\nu])-1$
 linearly independent solutions over the field of real numbers.

2) For   $n=[\nu]$  for the solvability of the  $(R-H)_0$  problem
it is necessary and sufficient that the single condition  (34) is
satisfied.

3)For  $n<[\nu]$ for the solvability of  the   $(R-H)_0$ problem it is
 necessary and sufficient that  $2n-[\nu]+1$  conditions of the  type (38)
 are satisfied, which are written with respect to the equation  (42).}

\vskip3 mm
{\bf References}

[1] Vekua I.N. Generalized analytic functions, Pergamon Press, London, Addison-Welsey, Reading, MA, 1962, MR 21\#7288; 27\#321.

[2] Usmanov Z.D. Infinitesimal bendings of surfaces of positive curvature with a flattening point.  Differential Geometry.
Banach Center Publications.  Warsaw. 1984. V. 12.  241--272.

[3] Usmanov Z.D. Infinitesimal bendings of surfaces of positive curvature with an isolated flattening point. Math. Sb. 1970,
 83(125): 4 (12),  596--615. (Russian).

[4] Mikhailov L.G. A New class of  singular integral equations and its applications to differential equations with a
singular coefficients. Noord Hoff Publishing House, Netherlands, Groningen, 1970.

[5] Usmanov Z.D. Generalized  Cauchy-Riemann  systems with a singular point. Sib.Math. 1973. J.14,  \textnumero 5 1076--1087. (Russian).

[6] Tungatarov A. On the theory of the Carleman-Vekua equation with a singular point. Russian Acad. Sci. Sb. Math. Vol.78 (1994), No.2. 357-365.

[7] Tungatarov A. On continuous solutions of of the Carleman-Vekua equation with a singular point. Soviet Math. Dokl. Vol.44 (1992). No.1. 175-178.

[8] Bizadze A.V. Basis of analytic functions theory of complex
variables  Nauka, Moskow, 1984. (Russian).

\end{document}